\documentclass[runningheads,envcountsame]{llncs}

\usepackage {amsmath, amssymb}
\usepackage[ruled,vlined]{algorithm2e}
\usepackage [top=1.5in, right=1in, bottom=1.5in, left=1in] {geometry}
\usepackage {pgfplots}
\pgfplotsset {compat=1.13}
\usepackage {tikz}
\usepackage{todonotes}
\usepackage{cite}
\usepackage{graphicx}
\usepackage {hyperref}
\hypersetup {
    colorlinks = true,
    linktoc = all,
    linkcolor = blue,
    citecolor = red
}

\usepackage {libertine}
\usepackage [T1]{fontenc}

\usepackage{authblk}

\numberwithin{theorem}{section}
\numberwithin{lemma}{section}
\numberwithin{proposition}{section}
\numberwithin{corollary}{section}
\numberwithin{example}{section}
\numberwithin{remark}{section}
\numberwithin{definition}{section}

\newcommand{\E}{\mathbb{E}}

\newcommand{\Z}{\mathbb{Z}}

\newcommand{\bS}{\mathbf{S}}
\newcommand{\bC}{\mathbf{C}}

\newcommand{\nauty}{\textsc{Nauty}}

\begin {document}

\title{Graph-like Scheduling Problems and Property B\thanks{With partial support from NSF Grant DMS-2039316.}}
\titlerunning{Graph-like scheduling}
\author{John Machacek\inst{1}}
\authorrunning{J. Machacek}
\institute{University of Oregon, Eugene OR  97403}
\maketitle

\begin {abstract}
Breuer and Klivans defined a diverse class of scheduling problems in terms of Boolean formulas with atomic clauses that are inequalities.
We consider what we call graph-like scheduling problems.
These are Boolean formulas that are conjunctions of disjunctions of atomic clauses $(x_i \neq x_j)$.
These problems generalize proper coloring in graphs and hypergraphs.
We focus on the existence of a solution with all $x_i$ taking the value of $0$ or $1$ (i.e. problems analogous to the bipartite case).
When a graph-like scheduling problem has such a solution, we say it has property B just as is done for $2$-colorable hypergraphs.
We define the notion of a $\lambda$-uniform graph-like scheduling problem for any integer partition $\lambda$.
Some bounds are attained for the size of the smallest $\lambda$-uniform graph-like scheduling problems without property B.
We make use of both random and constructive methods to obtain bounds.
Just as in the case of hypergraphs finding tight bounds remains an open problem.
\end {abstract}

\section {Introduction}
Various formulations of optimal job scheduling are frequently studied problems in computer science, operations research, and combinatorial optimization (see e.g.~\cite{GrahamSurvey, Graham}).
We will work with a formulation from Breuer and Klivans~\cite{BK16} which defines a \emph{scheduling problem} on $n$ elements to be a Boolean formula $\bS$ over atomic formulas $(x_i \leq x_j)$ for $1 \leq i,j \leq n$.
This means we can take arbitrary combinations of $(x_i \leq x_j)$ for $1 \leq i,j \leq n$ using \emph{and} denoted $\wedge$, \emph{or} denoted $\vee$, along with \emph{no}t denoted $\neg$.
A \emph{solution} to a scheduling problem is an assignment of a nonnegative integer to each $x_i$ that makes the Boolean formula $\bS$ true.
We can think of each $x_i$ as a task to be performed, the Boolean formula $\bS$ as constraints, and the assignments of nonnegative integers as a scheduling of the tasks into discrete time slots.
These problems contain several types of scheduling in the literature and have connections to algebra, combinatorics, and discrete geometry.
\\

This Boolean formula model of scheduling overlaps with other versions of scheduling in the literature of which we now name a few.
Perhaps, most notably it contains scheduling with precedence constraints\footnote{In our setup we would have access to as many processing machines as needed and all jobs would take the same unit time.} (see e.g.~\cite[4-3]{theoryofscheduling}).
An input to such a problem is \emph{precedence graph} which is a directed acyclic graph (DAG) where the vertices of the DAG are the jobs and if $(i,j)$ is an arc in the DAG then job $i$ must be completed before job $j$ begins.
For example, the DAG with arcs $\{(1,2), (1, 3), (3,4)\}$ would enocde as $(x_1 < x_2) \wedge (x_1 < x_3) \wedge (x_3 < x_4)$ in our notation.
Note that since $(x_i < x_j) = \neg (x_j \leq x_i)$ we can use strict inequality in our scheduling problems.
Scheduling with and/or precedence constraints has been considered where a job cannot start before at least one other job from a certain group~\cite{andorprec}.
Here the input is DAG where vertices are pairs $(S,j)$ with $S$ a set of jobs, and at least one job from $S$ must take place before job $j$.
For example, the input $\{(\{1,2\}, 3), (\{3,4\}, 1)\}$ becomes $((x_1 < x_3) \vee (x_2 < x_3)) \wedge ((x_3 < x_1) \vee (x_4 < x_1))$.
More recently, so-called \emph{s-precedence constraints} have been studied where these constraints require a job $j$ not to begin any earlier than some other job $i$ begins~\cite{sprec}.
Such constraints have applications in first-come, first-serve queuing and can be modeled using weak inequality in our setting.\\

We now move to consider scheduling problems that are similar to graph and hypergraph coloring.
It has long been known that there is a connection between graph coloring and the timetabling problem which asks for a way to schedule and avoid conflicts (see e.g.~\cite{Welsh,Wood}).
Let $\Z_{\geq 0}$ denote the nonnegative integers, then a \emph{coloring} of a hypergraph $H = (V,E)$ is just a map $\phi: V \to \Z_{\geq 0}$.
Here $V$ is a set of \emph{vertices} and $E$ is a set of \emph{hyperedges} which are simply subsets of $V$.
The case of a usual graph is when all subsets in $E$ have size $2$.
A \emph{proper coloring} is a coloring $\phi: V \to \Z_{\geq 0}$ such that $|\phi^{-1}(A)| > 1$ for all $A \in E$.
We refer to a nonnegative integer as a \emph{color} and  call a hyperedge $A$ \emph{monochromatic} with respect to a given coloring $\phi$ when $| \phi^{-1}(A)| = 1$.
So, a proper coloring is an assignment of colors to the vertices of the hypergraph such that no hyperedge is monochromatic.
If $H$ has a proper coloring such that $|\phi(V)| = k$, then $H$ is said to be \emph{$k$-colorable}.\\

In~\cite{Pluri} \emph{graph-like scheduling problems} on $n$ elements were defined to be Boolean formulas which are a conjunction of disjunctions over atomic formulas $(x_i \neq x_j)$ for $1 \leq i, j, \leq n$.
Since  
\[(x_i \neq x_j) = (x_i > x_j) \vee (x_j > x_i)\]
we see that graph-like scheduling problems are a special case of scheduling problems.
For example, the scheduling problem
\[(x_1 \neq x_2) \wedge ((x_1 \neq x_3) \vee (x_3 \neq x_4)) \wedge ((x_1 \neq x_4) \vee (x_2 \neq x_3))\]
is a graph-like scheduling problem.
One solution to this scheduling problem is to set $x_1 = 0$ and $x_i = 1$ for $1 < i \leq 4$.\\

Graph-like scheduling problems include hypergraph coloring as we now explain.
Define $[n] := \{1,2,\dots, n\}$ for any positive integer $n$.
For any subset $A \subseteq [n]$ define $A(i)$ to be the $i$th smallest element of $A$ so that $A = \{A(1) < A(2) < \cdots < A(|A|)\}$.
For any hypergraph $H = ([n], E)$ we define the graph-like scheduling problem\footnote{There are other logically equivalent ways to write down this scheduling problem. Using our notation $A(i)$ allows a canonical way to write it using a minimal number of clauses. For example, one could also use $(A(i) \neq A(j))$ for all $1 \leq i < j \leq |A|$.}
\[\bS_H := \bigwedge_{A \in E} \bigvee_{i = 1}^{|A|-1} (x_{A(i)} \neq x_{A(i+1)})\]
which models properly coloring $H$.
Indeed the clause 
\[ \bigvee_{i = 1}^{|A|-1} (x_{A(i)} \neq x_{A(i+1)})\]
is false only when $x_{i} = x_{j}$ for all $i,j \in A$ (i.e., when $A$ is monochromatic).
Hence, letting $x_i$ be assigned $\phi(i)$, the color received by vertex $i \in [n]$, the graph-like scheduling problem $\bS_H$ is true exactly when $\phi$ is a proper coloring of $H$.
We will focus of generalizing hypergraph coloring, but in~\cite[Section 5]{Pluri} it is shown graph-like scheduling problems also generalize other coloring problems such as orientated coloring~\cite{oriented} and star coloring~\cite{acyclic}.\\

This definition of scheduling problems in terms of Boolean formulas has been treated with algebraic and geometric combinatorics in~\cite{BK16} where associated quasisymmetric functions, simplicial complexes, and Ehrhart theory were considered. 
Further work on these scheduling problems in an algebraic direction has been done using Hopf algebras and Hopf monoids~\cite{aval, mario, white}.
Our present consideration is an extremal problem.
In particular, we wish to determine when the size of a graph-like scheduling problem requires it to have a solution with each $x_i$ being assigned $0$ or $1$.
This can be thought of as determining when we must have certain terms in an expansion of the quasisymmetric function mentioned above.
Alternatively, we can think of this as the existence of a solution in the hypercube with $0$-$1$ vector vertices from the Ehrhart theory perspective.
However, here we will not consider directly these algebraic or geometric notions.
Our results and techniques will be combinatorial and probabilistic.\\

Determining if a graph is $2$-colorable (i.e. bipartite) is easy and can be solved in linear time.
In contrast, testing if a hypergraph is $2$-colorable is difficult, and even for hypergraphs such that each hyperedge has size $3$ this is an NP-complete problem~\cite{L73}.
When phrased in terms of families of sets testing for $2$-colorability is known as \emph{the set splitting problem} which is one of Garey and Johnson's NP-complete problems~\cite{GJ}.
A hypergraph that is $2$-colorable is said to have \emph{property B}.
This is a long-studied property of hypergraphs which was named by Miller~\cite{Miller} with the B being for Bernstein~\cite{B} who first considered this property.
We say a graph-like scheduling problem has \emph{property B} if it has a solution so that each $x_i$ is assigned $0$ or $1$.
Of course, we can also make this definition for a general scheduling problem, but we will focus only on graph-like scheduling problems.\\

If $H = (V,E)$ is a hypergraph such that $|A| = k$ for all $A \in E$, then $H$ is called \emph{$k$-uniform}.
Let $m(k)$ denote the minimum number of hyperedges in a $k$-uniform hypergraph that does not have property B.
The problem of computing $m(k)$ was posed by Erd\H{o}s and Hajnal~\cite{EH61}.
Erd\H{o}s~\cite{combprobII} gave the best known general upper bound of $m(k) < k^2 2^k$ while the best known general lower bound is $m(k) = \Omega\left(\sqrt{\frac{k}{\log(k)}} 2^k\right)$~\cite{RS2000, CK15}.
There has also been work on finding bounds for particular small values of $k$~\cite{AASS2020}, and a survey of results as well as the history of the problem can be found in~\cite{survey, survey2}.
The only known exact values are $m(1) = 1$, $m(2) = 3$, $m(3) = 7$, and $m(4) = 23$.
The $k$-uniform hypergraphs without property B realizing the values of $m(k)$ for the cases of $1 \leq k \leq 3$ are a single vertex, the complete graph $K_3$, and the Fano plane.
For showing $m(4) = 23$ we have isomorphic examples from Seymour~\cite{S74} and Toft~\cite{T75} combined with an exhaustive computer search by \"{O}sterg\aa rd~\cite{O2014}.\\

In the next section, we define $\lambda$-uniform graph-like scheduling problems in a way that generalizes coloring in $k$-uniform hypergraphs where here $\lambda$ is an integer partition.
We then initiate the study of the analog of $m(k)$ for $\lambda$-uniform graph-like scheduling problems which we denote\footnote{We add the $s$ for ``scheduling'' to not risk confusion the notation with $m(k)$.} by $ms(\lambda)$.
In Theorem~\ref{thm:lower} we are able to give a lower bound for $ms(\lambda)$ using the probabilistic method.
We compute $ms((2,2))$ and then in Theorem~\ref{thm:2c} we give an upper bound for $ms((2c,2c))$ in terms of $ms((2,2))$ and $m(c)$.
In the final section, we define and investigate $k$-regular $\lambda$-uniform graph-like scheduling problems using the Lov\'{a}sz local lemma.

\section{$\lambda$-uniform Graph-like Scheduling Problems}

For a disjunction 
\[\bC = \bigvee_{(i,j) \in I} (x_i \neq x_j)\]
we define a graph $G_{\bC}$ with vertex set 
\[V_{\bC} = \{i : (i,j) \in I \text{ or } (j,i) \in I \text{ for some } j \}\]
and edge set
\[E_{\bC} = \{ij : (i,j) \in I\}\]
where we assume $(i,i) \not\in I$ for any $i \in [n]$ since $(x_i \neq x _i)$ is equivalent to false.
Next define $|\bC| := |V_{\bC}|$ and $\lambda(\bC) \vdash |\bC|$ to be the integer partition given by the sizes of the connected components of $G_{\bC}$.
Recall, an \emph{integer partition} $\lambda \vdash k$ is a nondecreasing sequence of positive integers summing to $k$.
For example, $(3,2,2,1) \vdash 8$ is an integer partition since $3+2+2+1 = 8$.
Notice a solution to $\bC$ is exactly an assignment such that for some component of $G_{\bC}$ there exists $i$ and $j$ in that component for which $x_i$ and $x_j$ are assigned different values.
That is, not all components of $G_{\bC}$ are monochromatic.
We identify clauses $\bC$ and $\bC'$ if $G_{\bC}$ and $G_{\bC'}$ have the same connected components.
This identification is logical equivalence in that if $G_{\bC}$ and $G_{\bC'}$ have the same connected components, then $\bC$ and $\bC'$ have the same truth value for any fixed assignment.
Given a graph-like scheduling problem
\[\bS = \bigwedge_{\alpha \in J} \bC_{\alpha}\]
if there is an integer $k$ and an integer partition $\lambda \vdash k$ such $\lambda(\bC_{\alpha}) = \lambda$ for each $\alpha \in J$ we say that $\bS$ is \emph{$\lambda$-uniform}.
We only will work with integer partitions that have no part equal to $1$ and use the notation $\lambda \vdash_{> 1} k$ to denote such partitions.\\

\begin{remark}
There is no loss of generality restricting to partitions with no part equal to $1$.
If $\lambda = (\lambda_1, \lambda_2, \dots, \lambda_{\ell-1}, 1)$, then each clause would need to have $(x_i \neq x_i)$ for some $i$ which can never be satisfied.
Hence, we could delete each instance of $(x_i \neq x_i)$ and look at the $\mu$-uniform graph-like scheduling problem for $\mu =  (\lambda_1, \lambda_2, \dots, \lambda_{\ell-1})$.
\end{remark}

Recall, a graph-like scheduling problem on $n$ elements $\bS$ has property B if there exists a solution to $\bS$ such that $x_i$ is assigned $0$ or $1$ for all $i \in [n]$.
Note this definition is consistent with hypergraph coloring in that a hypergraph $H$ has property B if and only if the associated graph-like scheduling problem $\bS_H$ has property B.
For a graph-like scheduling problem
\[\bS = \bigwedge_{\alpha \in J} \bC_{\alpha}\]
we say the \emph{size} of $\bS$ is $|J|$ and write $|\bS| = |J|$.

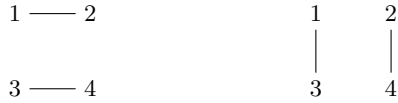
\begin{figure}
    \centering
    \begin{tikzpicture}
    \node (a) at (-2,1) {$1$};
    \node (b) at (-1,1) {$2$};
    \node (c) at (-2,0) {$3$};
    \node (d) at (-1,0) {$4$};
    \draw (a) to (b);
    \draw (c) to (d);

    \node (a) at (2,1) {$1$};
    \node (b) at (3,1) {$2$};
    \node (c) at (2,0) {$3$};
    \node (d) at (3,0) {$4$};
    \draw (a) to (c);
    \draw (b) to (d);
    \end{tikzpicture}
    \caption{The graphs $G_{\bC_1}$ on the left and $G_{\bC_2}$ on the right for $\bC_ 1 (x_1 \neq x_2) \vee (x_3 \neq x_4)$ and  $\bC_2 = (x_1 \neq x_3) \vee (x_2 \neq x_4)$.}
    \label{fig:graph}
\end{figure}

\begin{example}[A $(2,2)$-uniform graph-like scheduling problem]
Let $\bS = \bC_1 \wedge \bC_2$ for $\bC_ 1 = (x_1 \neq x_2) \vee (x_3 \neq x_4)$ and $\bC_2 = (x_1 \neq x_3) \vee (x_2 \neq x_4)$.
The graphs $G_{\bC_1}$ and $G_{\bC_2}$ are shown in Figure~\ref{fig:graph}.
The graph-like scheduling problem $\bS$ is $(2,2)$-uniform and $|\bS| = 2$.
Here $x_1 = 0$ and $x_2 = x_3 = x_4 = 1$ is a solution exhibiting that this $(2,2)$-uniform graph-like scheduling  problem has property B.
\end{example}

For any $\lambda \vdash_{> 1} k$ let $ms(\lambda)$ denote the minimum possible size of a $\lambda$-uniform graph-like scheduling problem that does not have property B.
Now the main problem we will focus on is attempting to compute, or at least bound, the values $ms(\lambda)$.
First, let us observe that $ms(\lambda)$ is a generalization of $m(k)$.

\begin{proposition}
For any $k \geq 2$ we have $m(k) = ms((k))$.
\label{prop:k}
\end{proposition}

\begin{proof}
For any $k$-uniform hypergraph $H$ we have that $\bS_H$ is a $(k)$-uniform graph-like scheduling problem.
Given a $(k)$-uniform graph-like scheduling problem
\[\bS = \bigwedge_{\alpha \in J} \bC_{\alpha}\]
the truth value of each clause $\bC_{\alpha}$ only depends on the connected components of $G_{\bC_{\alpha}}$.
By the $(k)$-uniform assumption $G_{\bC_{\alpha}}$ is connected and has $k$ vertices.
Hence, it follows $\bS$ is equivalent to $\bS_H$ were $H$ has hyperedge set $\{V_{\bC_{\alpha}} : \alpha \in J \}$.
Therefore, the proposition is proven.
\qed \end{proof}

Given two integer partitions $\lambda, \mu \vdash k$ we say that $\lambda$ \emph{refines} $\mu$ if $\mu$ can possibly be obtained by iterating the process of taking two parts of $\lambda$ and adding them to get one new larger part and shorter integer partition.
For example, $(4,3,2,2)$ refines $(5,4,2)$ and $(6,5)$.
We find the quantity $ms(\lambda)$ respects ordering by refinement.

\begin{proposition}
If $\lambda, \mu \vdash_{> 1} k$ and $\lambda$ refines $\mu$, then $ms(\lambda) \leq ms(\mu)$.
In particular, $ms(\lambda) \leq m(k)$ for all $\lambda \vdash_{> 1} k$.
\label{prop:refines}
\end{proposition}

\begin{proof}
The second statement follows from the first using Proposition~\ref{prop:k} and the fact that $(k)$ is refined by all integer partitions of $k$.
It then suffices to prove the first statement in the case that $\mu$ is obtained adding two parts of $\lambda$ since any relation by refinement comes from iterating the joining of two parts.
Assume $\mu$ comes from $\lambda$ by adding parts of size $a$ and $b$.
Take a $\lambda$-uniform graph-like scheduling problem
\[\bS = \bigwedge_{\alpha \in J} \bC_{\alpha}\]
and consider the $\mu$-uniform graph-like scheduling problem
\[\bS' = \bigwedge_{\alpha \in J} \bC'_{\alpha}\]
where $\bC'_{\alpha} = \bC_{\alpha} \vee (x_{i_{\alpha}} \neq x_{j_{\alpha}})$ such that $i_{\alpha}$ is in a component of $G_{\bC_{\alpha}}$ of size $a$ and $j_{\alpha}$ is in a component of $G_{\bC_{\alpha}}$ of size $b$ distinct from the component containing $i_{\alpha}$.
For a given assignment we have that if $\bS'$ is false, then $\bS$ is false.
So, if $\bS'$ does not have property $B$, then $\bS$ does not have property $B$.
The proposition then follows.
\qed \end{proof}

\subsection{A Lower Bound for $ms(\lambda)$}

We now give a simple lower bound that comes from a random assignment that agrees with the lower bound of Erd\H{o}s~\cite{E63} for partitions with length equal to $1$ (i.e. for $k$-uniform hypergraphs).
It would be desirable to find a bound using more sophisticated techniques.
Possibly something similar to what was done to improve the bound for $m(k)$ using a random ordering of the vertices in~\cite{RS2000, CK15} could be used.
However, the bound in Theorem~\ref{thm:lower} states that $4 \leq ms((2,2))$, and we will find that $ms((2,2)) = 5$ in Lemma~\ref{lem:22}.

\begin{theorem}
If\, $\bS$ is a $\lambda$-uniform graph-like scheduling problem with size $|\bS| < 2^{k - \ell}$ where $\lambda \vdash_{> 1} k$ is partition of length $\ell$, then $\bS$ has property B.
Hence, $2^{k - \ell} \leq ms(\lambda)$ for any $\lambda \vdash_{> 1} k$ of length $\ell$.
\label{thm:lower}
\end{theorem}

\begin{proof}
Let
\[\bS = \bigwedge_{\alpha \in J} \bC_{\alpha}\]
be our $\lambda$-uniform graph-like scheduling problem on $n$ elements.
Consider assigning $0$ or $1$ to each $x_i$ uniformly at random so that
\[\Pr[x_i \text{ is assigned } 0] = \frac{1}{2} = \Pr[x_i \text{ is assigned } 1]\]
for all $i \in [n]$.
For each $\alpha \in J$ let $Y_{\alpha}$ be the random variable such that $Y_{\alpha} = 1$ if $\bC_{\alpha}$ is false and $Y_{\alpha} = 0$ if $\bC_{\alpha}$ is true.
Now we look at an expected value
\begin{align*}
\E \Big[ \big| \{\bC_{\alpha} : \alpha \in J, \bC_{\alpha} \text{ is false}\} \big| \Big] &= \E \left[ \sum_{\alpha \in J} Y_{\alpha} \right]\\
&= \sum_{\alpha \in J} \E [ Y_{\alpha} ]\\
&= \sum_{\alpha \in J} \frac{2^{\ell}}{2^k}\\
&= \frac{|J|}{2^{k-\ell}}\\
&< 1
\end{align*}
which holds since $|J| = |\bS| < 2^{k - \ell}$.
Now since the expected value is less than $1$ and each $Y_{\alpha}$ only takes values $0$ or $1$ it follows there exists some assignment so that $Y_{\alpha} = 0$ for all $\alpha \in J$.
Such an assignment makes $\bS$ true and therefore $\bS$ has property $B$.
\qed \end{proof}

\subsection{An Upper Bound for $ms((2c,2c))$}

We first compute the smallest case which is not a uniform hypergraph problem.
Note that by Proposition~\ref{prop:refines} and the fact that $m(4) = 23$ we have that $ms((2,2)) \leq 23$.
However, it turns out $ms((2,2))$ is much smaller and close to bound in Theorem~\ref{thm:lower}.
We will use exhaustive computation to help compute the value of $ms((2,2))$, but first we prove some results to help reduce the computation need.\\

We define a \emph{path} to be a sequence of distinct vertices and hyperedges $v_0, e_1, v_1, e_2, \dots, v_{\ell-1}, e_{\ell}, v_\ell$ so that $v_i \in V$ for $0 \leq i \leq \ell$, $e_i \in E$ for $1 \leq i \leq \ell$, $v_0 \in e_0$, $v_{\ell} \in e_{\ell}$, and $v_i \in e_i \cap e_{i+1}$ for $0 < i < \ell$.
A hypergraph is \emph{connected} if there is a path between any two vertices.
A path is a \emph{cycle} if $\ell > 1$ and $v_0 = v_{\ell}$.
A hypergraph without any cycle is called \emph{acylic}.
The lemma we will need uses the notion of a \emph{hypertree} which is a connected acyclic hypergraph.

\begin{example}[Cycles in hypergraphs]
    Consider the hypergraph $H = (V,E)$ where $V = \{1,2,3,4,5\}$ and 
    \[E = \{\{1,2,3\}, \{1,5\} \{2,3,4\}, \{4,5\}\}.\]
    This hypergraph has the cycle $2, \{1,2,3\}, 3, \{2,3,4\}, 2$ as well as the cycle 
    \[1, \{1,5\}, 5, \{4,5\}, 4 \{2,3,4\}, 2, \{1,2,3\}, 1.\]
\end{example}

The example shows that if we ever have two hyperedges $e \neq e'$ with $|e \cap e'| \geq 2$ then we have a cycle.
Indeed, if $v, v' \in e \cap e'$ with $v \neq v'$, then $v, e, v', e', v$ is a cycle.
This means any hypertree must have $|e \cap e'| \leq 1$ for any $e \neq e'$.
A hypergraph satisfying this condition on hyperedge intersection size is known as a \emph{linear hypergraph}.
Any graph-like scheduling problem $\bigwedge_{\alpha \in J} \bC_{\alpha}$ determines a hypergraph $([n], \{V_\alpha : \alpha \in J\})$. 
We will need to consider all $\lambda$-uniform graph-like scheduling problems giving rise to a fixed uniform hypergraph.
This can be done considering all ways of partitioning each hyperedge to form a clause as the next example shows.

\begin{example}[Moving between hypergraphs to scheduling problems]
    Consider the $4$-uniform hypergraph $H = (V,E)$ where $V = \{1,2,3,4,5,6\}$ and $E = \{\{1,2,3,4\}, \{3,4,5,6\}\}$.
    To find all $(2,2)$-uniform graph-like scheduling problems corresponding to this hypergraph we must partition each hyperedge into $2$ disjoint subsets of size $2$.
    For the hyperedge $\{1,2,3,4\}$, since the order of the blocks does not matter, we have the $3$ such partitions
    \[\{1,2\}/\{3,4\}, \{1,3\}/\{2,4\}, \{1,4\}/\{2,3\}\]
    and similarly for the hyperedge $\{3,4,5,6\}$.
    Thus we get the following possible $(2,2)$-uniform graph-like scheduling problems.
    \begin{align*}
        ((x_1 \neq x_2) \vee (x_3 \neq x_4)) \wedge ((x_3 \neq x_4) \vee (x_5 \neq x_6))\\
        ((x_1 \neq x_3) \vee (x_2 \neq x_4)) \wedge ((x_3 \neq x_4) \vee (x_5 \neq x_6))\\
        ((x_1 \neq x_4) \vee (x_2 \neq x_3)) \wedge ((x_3 \neq x_4) \vee (x_5 \neq x_6))\\
        ((x_1 \neq x_2) \vee (x_3 \neq x_4)) \wedge ((x_3 \neq x_5) \vee (x_4 \neq x_6))\\
        ((x_1 \neq x_3) \vee (x_2 \neq x_4)) \wedge ((x_3 \neq x_5) \vee (x_4 \neq x_6))\\
        ((x_1 \neq x_4) \vee (x_2 \neq x_3)) \wedge ((x_3 \neq x_5) \vee (x_4 \neq x_6))\\
        ((x_1 \neq x_2) \vee (x_3 \neq x_4)) \wedge ((x_3 \neq x_6) \vee (x_4 \neq x_5))\\
        ((x_1 \neq x_3) \vee (x_2 \neq x_4)) \wedge ((x_3 \neq x_6) \vee (x_4 \neq x_5))\\
        ((x_1 \neq x_4) \vee (x_2 \neq x_3)) \wedge ((x_3 \neq x_6) \vee (x_4 \neq x_5))\\       
    \end{align*}
    as the $9$ graph-like scheduling problems corresponding to our fixed hypergraph.
    \label{ex:hypersched}
\end{example}

A $k$-uniform hypergraph on $n$ vertices with $m$ hyperedges is a hypertree if and only if it is connected and $km = n-1$~\cite[Lemma 4.2]{Pluri}.
This will be useful to determine when we have a hypertree.
We call a hyperedge a \emph{leaf} if all but at most one of its vertices have degree equal to $1$ in the hypergraph.
Note that this generalizes the notion of an edge that contains a leaf vertex of a graph.
It turns out that any hypertree necessarily has a leaf.

\begin{lemma}
    If $H$ is a $k$-uniform hypertree with at least one hyperedge, then $H$ has at least one leaf.
    \label{lem:leaf}
\end{lemma}
\begin{proof}
Assume $H$ has no leaf.
Then we will show we can find a cycle in $H$.
Choose any vertex $v_0$ and hyperedge $e_1$ with $v_0 \in e_1$.
Now for any $i \leq 1$ since $e_i$ is not a leaf we have $v_i \in e_i$ with degree at least $2$ where $v_i \neq v_{i-1}$.
In this way, we can build a walk $v_0, e_1, v_1, e_2, v_2, \dots$ which will eventually contain a cycle.
\qed \end{proof}

\begin{lemma}
    If a $\lambda$-uniform graph-like scheduling problem $\bS$ corresponds to a hypertree, then $\bS$ has property B.
    \label{lem:tree}
\end{lemma}
\begin{proof}
    The lemma follows by induction on the size of the hypertree by removing a leaf that exists by Lemma~\ref{lem:leaf}.
    If $\bS$ corresponds to a hypertree with a single hyperedge, then $\bS$ has property $B$.
    If $\bS$ is a $\lambda$-uniform graph-like scheduling problem corresponding to a hypertree with more than one hyperedge consider a leaf of the hypertree.
    Remove the clause corresponding to this leaf along with all degree $1$ vertices in the leaf to obtain a smaller $\lambda$-uniform graph-like scheduling problem $\bS$' that has property B by induction.
    Take a solution to $\bS$ using only values $0$ and $1$.
    Now because the clause removed corresponded to a leaf only a single $x_i$ in the clause is assigned by the solution to $\bS$'.
    Hence, this solution can be extended to $\bS$, and therefore $\bS$ has property B.
\qed \end{proof}

\begin{lemma}
It is the case that $ms((2,2)) = 5$.
\label{lem:22}
\end{lemma}
\begin{proof}
By Theorem~\ref{thm:lower} we know that $ms(2,2) \geq 4$.
We may argue that $ms((2,2)) > 4$ by direct computation.
To do this we can first, for some $n$, generate connected $4$-uniform hypergraphs on $n$ vertices with $4$ hyperedges using \nauty~\cite{nauty}.
Then for each hyperedge $\{i, j, k, \ell\}$ we enumerate each set partition with $2$ parts of size $2$.
So, the set partition $\{i,j\} / \{k, \ell\}$ would correspond to the clause $(x_i \neq x_j) \vee (x_k \neq x_{\ell})$ as in Example~\ref{ex:hypersched}.
For each such hypergraph along with set partitions for each hyperedge, we may then check the corresponding graph-like scheduling problem has property B by testing possible assignments with $0$ or $1$ to each variable.
Note that we will have to repeat for multiple values of $n$, but we can stop after $n = 12$. 
A $4$-uniform hypergraph with $n > 13$ will not be connected.
We may restrict our attention to the connected case since a graph-like scheduling problem will have property B if and only if each connected component does.
When $n=13$, any connected hypergraph will be a hypertree since $m(k-1) = n-1$ for $m=4$, $k=4$, and $n = 13$.
By Lemma~\ref{lem:tree} all such graph-like scheduling problems will have property B when $n=13$.\\

Then to establish that $ms((2,2)) = 5$ we may check that
\begin{align*}
&((x_1 \neq x_2) \vee (x_3 \neq x_4)) \wedge ((x_1 \neq x_3) \vee (x_2 \neq x_5)) \wedge ((x_1 \neq x_5) \vee (x_2 \neq x_4))\\
 &\quad\quad \wedge ((x_1 \neq x_4) \vee(x_3 \neq x_5)) \wedge ((x_2 \neq x_3) \vee (x_4 \neq x_5))
\end{align*}
does not have property B.
This can be done by simply testing all possible assignments for the variables with $0$ or $1$.
\qed \end{proof}

For an integer partition $\lambda = (\lambda_1, \lambda_2, \dots, \lambda_{\ell}) \vdash k$ for positive integer $c$  let $c \cdot \lambda = (c\lambda_1, c\lambda_2, \dots, c\lambda_{\ell}) \vdash ck$.
We now give the following lemma which generalizes a construction of Abbott and Moser~\cite{AbbottMoser} from hypergraphs to graph-like scheduling problems.

\begin{lemma}
For any $\lambda \vdash_{> 1} k$ and $c > 0$ we have $ms(c \cdot \lambda) \leq ms(\lambda)m(c)^{k}$.
\label{lem:construction}
\end{lemma}
\begin{proof}
Let
\[\bS = \bigwedge_{\alpha \in J} \bC_{\alpha}\]
be a $\lambda$-uniform graph-like scheduling problem on $n$ elements without property B.
For each $i \in [n]$ let $H_i$ be a $c$-uniform hypergraph without property B such that when $i \neq j$ the hypergraphs $H_i$ and $H_j$ are disjoint.
We may assume $\bS$ and $H_i$ for $i \in [n]$ have size $ms(\lambda)$ and $m(c)$ respectively.\\

We construct the $c\cdot \lambda$-uniform graph-like scheduling problem on $N$ elements
\[\bS' = \bigwedge_{\alpha' \in J'} \bC'_{\alpha'}\]
where $\alpha' \in J'$ is the data of a choice $\alpha \in J$ along with a choice $e_i \in H_i$ for each $i \in V_{\bC_{\alpha}}$, and the clause $\bC'_{\alpha'}$ is determined by (the components of) the graph $G_{\bC'_{\alpha'}}$ that has vertex set 
\[V_{\bC'_{\alpha'}} =  \bigcup_{i \in V_{\bC_{\alpha}} }e_i\]
and a component on vertices
\[\bigcup_{i \in A}e_i\]
for each $A$ which is a vertex set of a component of $G_{\bC_{\alpha}}$.
Here we can assume that $[N]$ equals the (disjoint) union of vertex sets of the hypergraphs $H_i$ over $i \in [n]$.\\

We claim that $\bS'$ does not have property $B$.
Consider any $f:[N] \to \{0,1\}$, we must show $f$ is not a solution to $\bS'$.
The map $f$ restricts to a coloring of each $H_i$, and since $H_i$ does not have property $B$ there exists a monochromatic hyperedge in $H_i$ which we denote $\tilde{e}_i$.
Now let 
\[\tilde{\bS} = \bigwedge_{\tilde{\alpha} \in \tilde{J}} \bC'_{\tilde{\alpha}}\]
where $\tilde{J} \subseteq J'$ indexes those clauses gotten by choosing $\alpha \in J$ along with $\tilde{e}_i \in H_i$ for each $i \in V_{\bC_{\alpha}}$.
It is clear that $\bS' = \tilde{\bS} \wedge \bS''$ for some $\bS''$, and thus if $f$ is not a solution to $\tilde{\bS}$ it also is not a solution to $\bS'$.
The truth value of $\tilde{\bS}$ with assignment given by $f$ is the same as the truth value $\bS$ with assignment $x_i = g(i)$ given by $g: [n] \to \{0,1\}$ where $g(i) = f(j)$ for $j \in \tilde{e_i}$.
The map $g$ is well defined because $f$ is constant on each $\tilde{e}_i$.
Since $\bS$ does not have property B it must be the case that $\bS$ is false with this assignment.
It follows that $\tilde{\bS}$ and also $\bS'$ are false with the assignment under consideration.
As $f$ was arbitrary we can conclude that $\bS'$ does not have property B and the lemma is proven.
\qed \end{proof}

\begin{theorem}
If $c > 0$, then $ms((2c,2c)) \leq 5 m(c)^4$.
\label{thm:2c}
\end{theorem}

\begin{proof}
The theorem follows immediately from Lemma~\ref{lem:22} and Lemma~\ref{lem:construction}.
\qed \end{proof}

\begin{remark}
To get bounds for specific values of $c$ one can take Theorem~\ref{thm:2c} along with state-of-the-art bounds for $m(c)$.
For small values of $c$, bounds for $m(c)$ are tabulated in~\cite[Table 1]{AASS2020}.
\end{remark}

\section{$k$-regular Graph-like Scheduling Problems}
A graph-like scheduling problem 
\[\bS = \bigwedge_{\alpha \in J} \bC_{\alpha}\]
on $n$ elements is \emph{$k$-regular} provided for all $i \in [n]$ we have $\big| \{\bC_{\alpha} : i \in V_{\bC_{\alpha}} \} \big| = k$.
In this section, we will need a version of the \emph{Lov\'{a}sz local lemma}.
The version we will use is as follows (see~\cite[Corollary 5.1.2]{probmethod}).
Let $A_1$, $A_2$, $\dots$, $A_n$ be events in a probability space.
Assume $\Pr[A_i] \leq p$ for each $i$ and additionally, each $A_i$ is mutually independent for all but at most $d$ other $A_j$. If $ep(d+1) \leq 1$, then $\Pr\left[ \bigwedge_{i=1}^n \overline{A_i}\right] > 0$.

\begin{theorem}
Let 
\[\bS = \bigwedge_{\alpha \in J} \bC_{\alpha}\]
be a $\lambda$-uniform graph-like scheduling problem for $\lambda \vdash_{> 1} k$ of length $\ell$ such that
\[\Big| \{ \beta : \beta \in J, \beta \neq \alpha, V_{\bC_{\alpha}} \cap V_{\bC_{\beta}} \neq \emptyset\} \Big| < \frac{2^{k - \ell}}{e} - 1\]
for all $\alpha \in J$, then $\bS$ has property $B$.
\label{thm:prereg}
\end{theorem}

\begin{proof}
Take $\bS$ as in the statement of the theorem, and uniformly at random assign $0$ or $1$ to each $x_i$ so that
\[\Pr[x_i \text{ is assigned } 0] = \frac{1}{2} = \Pr[x_i \text{ is assigned } 1]\]
for all $i \in [n]$.
For each $\alpha \in J$ let $A_{\alpha}$ be the event that $\bC_{\alpha}$ is not satisfied, then
\[\Pr(A_{\alpha}) = 2^{\ell - k}\]
and let $p = 2^{\ell - k}$ be the common event probability.
The event $A_{\alpha}$ is independent from $A_{\beta}$ unless $V_{\bC_{\alpha}} \cap V_{\bC_{\beta}} \neq \emptyset$.
So, by the assumption in the theorem each $A_{\alpha}$ is independent from all but at most $\frac{2^{k - \ell}}{e} - 1$ of the $A_{\beta}$.
Noting that
\[e\cdot 2^{\ell - k} \cdot \left(\frac{2^{k - \ell}}{e}\right) = 1\]
and by applying the Lov\'asz local lemma it follows there is a nonzero probability that none of the events $A_{\alpha}$ for $\alpha \in J$ occur.
Therefore $\bS$ has property B.
\qed \end{proof}

\begin{table}
\centering
\begin{tabular}{|c|c|} \hline
$k$ & $e(k(k-1) + 1) 2^{-k}$ \\ \hline
8 & 0.60524244 \\ \hline
\textbf{9} & \textbf{0.38756753} \\ \hline
\textbf{10} & \textbf{0.24156606} \\ \hline
11 & 0.14732875 \\ \hline
\textbf{12} & \textbf{0.088264522} \\ \hline
\textbf{13} & \textbf{0.052095977} \\ \hline
\textbf{14} & \textbf{0.030361668} \\ \hline
15 & 0.017503585 \\ \hline
\textbf{16} & \textbf{0.0099961231} \\ \hline
\textbf{17} & \textbf{0.0056617046} \\ \hline
\textbf{18} & \textbf{.0031834126} \\ \hline
\textbf{19} & \textbf{0.0017783559} \\ \hline
20 & 0.00098768747 \\ \hline
\end{tabular}
\label{tbl:reg}
\caption{Values of the quantity from Corollary~\ref{cor:reg}. Bold rows indicate that the row corresponds to the smallest value of $k$ for which the quantity is less than some power of $2$.}
\end{table}

\begin{corollary}
If $\bS$ is a $\lambda$-uniform $k$-regular graph-like scheduling problem and
\[e(k(k-1) + 1) 2^{-k} < 2^{-\ell}\]
where $\lambda \vdash_{> 1} k$ has length $\ell$, then $\bS$ has property B.
Hence, for any fixed $\ell > 0$ there exists $K > 0$ such that for all $k \geq K$ any $\lambda$-uniform $k$-regular graph-like scheduling problem where $\lambda \vdash_{> 1} k$ is of length $\ell$ has property B.
\label{cor:reg}
\end{corollary}
\begin{proof}

Let $\bS$ be a $\lambda$-uniform $k$-regular graph-like scheduling problem for $\lambda \vdash_{> 1} k$ which has length $\ell$.
Write
\[\bS = \bigwedge_{\alpha \in J} \bC_{\alpha}\]
then for $\alpha \in J$ fixed
\[\Big| \{ \beta : \beta \in J, \beta \neq \alpha, V_{\bC_{\alpha}} \cap V_{\bC_{\beta}} \neq \emptyset\} \Big| \leq k(k-1)\]
since $\bS$ is $k$-regular and $\lambda$-uniform where $\lambda \vdash_{> 1} k$.
Applying Theorem~\ref{thm:prereg} and seeing that $k(k-1) < \frac{2^{k - \ell}}{e} - 1$ holds if and only if $e(k(k-1) + 1) 2^{-k} < 2^{-\ell}$ the corollary is proven.
\qed \end{proof}

The argument above is a variation of a known corresponding argument for hypergraphs (see e.g.~\cite[Section 5.2]{probmethod})
Values of $e(k(k-1) + 1) 2^{-k}$  are tabulated in Table~\ref{tbl:reg}.
These values can be used to provide a value of $K$ in Corollary~\ref{cor:reg}.
For example, when $\lambda = (\lambda_1, \lambda_2, \lambda_3) \vdash_{> 1} k$ from the table we find any $k$-regular $\lambda$-uniform graph-like scheduling problem has property B whenever $k \geq 12$.
This is because $k=12$ is the first row for which $e(k(k-1) + 1) 2^{-k} < 2^{-3} = 0.125$; so, we obtain $K = 12$ is the corollary.
However, these values of $K$ need not be optimal.
Looking at $\ell=1$ we obtain that all $k$-uniform $k$-regular hypergraphs have property B whenever $k \geq 9$, but it is known this can be improved to $k \geq 4$~\cite{Thom92}.

\bibliographystyle{splncs04}
\bibliography{refs.bib}

\begin{thebibliography}{10}
\providecommand{\url}[1]{\texttt{#1}}
\providecommand{\urlprefix}{URL }
\providecommand{\doi}[1]{https://doi.org/#1}

\bibitem{AbbottMoser}
Abbott, H.L., Moser, L.: On a combinatorial problem of {E}rdős and {H}ajnal.
  Canad. Math. Bull.  \textbf{7},  177--181 (1964).
  \doi{10.4153/CMB-1964-016-9}

\bibitem{AASS2020}
Aglave, S., Amarnath, V.A., Shannigrahi, S., Singh, S.: Improved bounds for
  uniform hypergraphs without property {B}. Australas. J. Combin.
  \textbf{76}(part 1),  73--86 (2020)

\bibitem{probmethod}
Alon, N., Spencer, J.H.: The probabilistic method. Wiley-Interscience Series in
  Discrete Mathematics and Optimization, Wiley-Interscience [John Wiley \&
  Sons], New York, second edn. (2000). \doi{10.1002/0471722154}, with an
  appendix on the life and work of Paul Erd\H{o}s

\bibitem{aval}
Aval, J.C., Bergeron, N., Machacek, J.: New invariants for permutations, orders
  and graphs. Adv. in Appl. Math.  \textbf{121},  102080, 30 (2020).
  \doi{10.1016/j.aam.2020.102080}

\bibitem{B}
Bernstein, F.: Zur theorie der trigonometrische reihen. Leipz. Ber.
  \textbf{60},  325--328 (1908)

\bibitem{BK16}
Breuer, F., Klivans, C.J.: Scheduling problems. J. Combin. Theory Ser. A
  \textbf{139},  59--79 (2016). \doi{10.1016/j.jcta.2015.11.001}

\bibitem{CK15}
Cherkashin, D.D., Kozik, J.: A note on random greedy coloring of uniform
  hypergraphs. Random Structures Algorithms  \textbf{47}(3),  407--413 (2015).
  \doi{10.1002/rsa.20556}

\bibitem{theoryofscheduling}
Conway, R., Maxwell, W.L., Miller, L.W.: Theory of scheduling. Addison-Wesley
  Pub. Co Reading, Mass (1967)

\bibitem{E63}
Erd\H{o}s, P.: On a combinatorial problem. Nordisk Mat. Tidskr.  \textbf{11},
  5--10, 40 (1963)

\bibitem{combprobII}
Erd\H{o}s, P.: On a combinatorial problem. {II}. Acta Math. Acad. Sci. Hungar.
  \textbf{15},  445--447 (1964). \doi{10.1007/BF01897152}

\bibitem{EH61}
Erd\H{o}s, P., Hajnal, A.: On a property of families of sets. Acta Math. Acad.
  Sci. Hungar.  \textbf{12},  87--123 (1961). \doi{10.1007/BF02066676}

\bibitem{GJ}
Garey, M.R., Johnson, D.S.: Computers and intractability. W. H. Freeman and
  Co., San Francisco, Calif. (1979), a guide to the theory of NP-completeness,
  A Series of Books in the Mathematical Sciences

\bibitem{andorprec}
Gillies, D.W., Liu, J.W.S.: Scheduling tasks with and/or precedence
  constraints. SIAM J. Comput.  \textbf{24}(4),  797--810 (1995).
  \doi{10.1137/S0097539791218664}

\bibitem{GrahamSurvey}
Graham, R.L., Lawler, E.L., Lenstra, J.K., Rinnooy~Kan, A.H.G.: Optimization
  and approximation in deterministic sequencing and scheduling: a survey. In:
  Interfaces between computer science and operations research ({P}roc.
  {S}ympos., {M}ath. {C}entrum, {A}msterdam, 1976). Math. Centre Tracts,
  vol.~99, pp. 169--214. Math. Centrum, Amsterdam (1978)

\bibitem{Graham}
Graham, R.L.: Combinatorial scheduling theory. In: Steen, L.A. (ed.)
  Mathematics Today Twelve Informal Essays, pp. 183--211. Springer New York,
  New York, NY (1978). \doi{10.1007/978-1-4613-9435-8\_8}

\bibitem{acyclic}
Gr{\"{u}}nbaum, B.: Acyclic colorings of planar graphs. Israel J. Math.
  \textbf{14},  390--408 (1973). \doi{10.1007/BF02764716}

\bibitem{sprec}
Kim, E.S., Posner, M.E.: Parallel machine scheduling with s-precedence
  constraints. IIE Transactions  \textbf{42}(7),  525--537 (2010).
  \doi{10.1080/07408171003670975}

\bibitem{L73}
Lov\'{a}sz, L.: Coverings and coloring of hypergraphs. In: Proceedings of the
  {F}ourth {S}outheastern {C}onference on {C}ombinatorics, {G}raph {T}heory,
  and {C}omputing ({F}lorida {A}tlantic {U}niv., {B}oca {R}aton, {F}la., 1973).
  pp. 3--12 (1973)

\bibitem{Pluri}
Machacek, J.: Plurigraph coloring and scheduling problems. Electron. J. Combin.
   \textbf{24}(2),  Paper No. 2.29 (2017). \doi{10.37236/6818}

\bibitem{nauty}
McKay, B.D., Piperno, A.: Practical graph isomorphism, {II}. J. Symbolic
  Comput.  \textbf{60},  94--112 (2014). \doi{10.1016/j.jsc.2013.09.003}

\bibitem{Miller}
{Miller}, E.W.: {On a property of families of sets}. {C. R. Soc. Sci. Varsovie,
  Cl. III}  \textbf{30},  31--38 (1937)

\bibitem{RS2000}
Radhakrishnan, J., Srinivasan, A.: Improved bounds and algorithms for
  hypergraph {$2$}-coloring. Random Structures Algorithms  \textbf{16}(1),
  4--32 (2000).
  \doi{10.1002/(SICI)1098-2418(200001)16:1<4::AID-RSA2>3.3.CO;2-U}

\bibitem{survey}
Raigorodskii, A.M., Cherkashin, D.D.: Extremal problems in hypergraph
  colouring. Uspekhi Mat. Nauk  \textbf{75}(1(451)),  95--154 (2020).
  \doi{10.4213/rm9905}

\bibitem{survey2}
Ra\u{\i}gorodski\u{\i}, A.M., Shabanov, D.A.: The {E}rd{\H{o}}s-{H}ajnal
  problem of hypergraph colorings, its generalizations, and related problems.
  Uspekhi Mat. Nauk  \textbf{66}(5(401)),  109--182 (2011).
  \doi{10.1070/RM2011v066n05ABEH004764}

\bibitem{O2014}
\"{O}sterg\aa rd, P.R.J.: On the minimum size of 4-uniform hypergraphs without
  property {$B$}. Discrete Appl. Math.  \textbf{163}(part 2),  199--204 (2014).
  \doi{10.1016/j.dam.2011.11.035}

\bibitem{mario}
Sanchez, M.: M\"{o}bius inversion as duality for {H}opf monoids. S\'{e}m.
  Lothar. Combin.  \textbf{82B},  Art. 91, 12 (2020)

\bibitem{S74}
Seymour, P.D.: A note on a combinatorial problem of {E}rd{\H{o}}s and {H}ajnal.
  J. London Math. Soc. (2)  \textbf{8},  681--682 (1974).
  \doi{10.1112/jlms/s2-8.4.681}

\bibitem{oriented}
Sopena, E.: Homomorphisms and colourings of oriented graphs: an updated survey.
  Discrete Math.  \textbf{339}(7),  1993--2005 (2016).
  \doi{10.1016/j.disc.2015.03.018}

\bibitem{Thom92}
Thomassen, C.: The even cycle problem for directed graphs. J. Amer. Math. Soc.
  \textbf{5}(2),  217--229 (1992). \doi{10.2307/2152767},
  \url{https://doi.org/10.2307/2152767}

\bibitem{T75}
Toft, B.: On colour-critical hypergraphs. In: Infinite and finite sets
  ({C}olloq., {K}eszthely, 1973; dedicated to {P}. {E}rd\H{o}s on his 60th
  birthday), {V}ol. {III}, pp. 1445--1457. Colloq. Math. Soc. J\'{a}nos Bolyai,
  Vol. 10 (1975)

\bibitem{Welsh}
Welsh, D.J.A., Powell, M.B.: {An upper bound for the chromatic number of a
  graph and its application to timetabling problems}. The Computer Journal
  \textbf{10}(1),  85--86 (1967). \doi{10.1093/comjnl/10.1.85}

\bibitem{white}
White, J.A.: Quasisymmetric functions from combinatorial {H}opf monoids and
  {E}hrhart theory. In: 28th {I}nternational {C}onference on {F}ormal {P}ower
  {S}eries and {A}lgebraic {C}ombinatorics ({FPSAC} 2016), pp. 1215--1226.
  Discrete Math. Theor. Comput. Sci. Proc., BC, Assoc. Discrete Math. Theor.
  Comput. Sci., Nancy ([2016] \copyright 2016)

\bibitem{Wood}
Wood, D.C.: {A technique for colouring a graph applicable to large scale
  timetabling problems}. The Computer Journal  \textbf{12}(4),  317--319
  (1969). \doi{10.1093/comjnl/12.4.317}

\end{thebibliography}

\end {document}